\documentclass[11pt]{amsart}
\usepackage[utf8]{inputenc}
\usepackage[english]{babel}
 \usepackage{float}
\usepackage{color}
\usepackage{graphicx}
\usepackage{graphics}
\usepackage{amsbsy, amsmath,amsfonts,amssymb,amsthm,amscd,amssymb,latexsym}

\numberwithin{equation}{section}

\usepackage{hyperref}
\hypersetup{pdfnewwindow=true,colorlinks=true,linkcolor=blue,citecolor=blue,
	filecolor=magenta,urlcolor=blue 
	}

\usepackage{graphicx,xcolor}			

\usepackage{amsthm}
\usepackage{multicol} 
\usepackage{graphicx}
\usepackage{fancyhdr}
\newtheorem{lema}{Lemma}[section]
\newtheorem{thm}{Theorem}[section]

\newtheorem{defi}{Definition}[section]
\newtheorem{pro}{Proposition}[section]
\newtheorem{remark}{Remark}[section]

\author{  Ronaldo Garcia   and Dimas Tejada }

\title[Principal lines on an ellipsoid]{  Principal Lines on an Ellipsoid in a Minkoswki Tridimensional Space} 
 
\date{ \today}

\begin{document}
\maketitle

\begin{abstract}
The description of principal lines of the ellipsoid on the 3-dimensional Minkowski space is established.  A global principal parametrization of a triple orthogonal system of quadrics is also achieved, and the focal set of the ellipsoid is sketched.
\end{abstract}


\section{Introduction}

The goal of this work is to describe the global behavior of principal lines of the ellipsoid in the three dimensional Minkwoswki space $\mathbb{R}^{2,1}$.  We recall that the concept of principal lines were introduced by G. Monge \cite{monge} and geometrically they can be characterized as the curves on the surface such that the ruled surface having the  rules being the  normal straight lines along the curve is a developable surface \cite[page 93]{Struik}.
 
\begin{figure}[H]
 		\centering	\includegraphics[width=0.7\textwidth,angle=0]{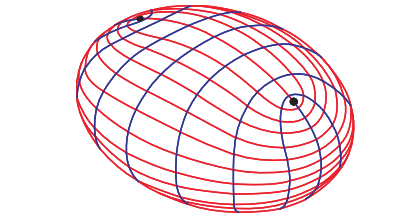}
 		\caption{ Principal lines on the triaxial ellipsoid. There are four umbilic points, the singularities. Also, there are four umbilic separatrices and other principal lines are closed.}	\label{fig:emonge}
 	\end{figure}

The principal lines of the ellipsoid with three  different axes  in the Euclidean space $\mathbb{R}^3$  are as illustrated in Fig. \ref{fig:emonge}.
In this case, the principal lines of the triaxial ellipsoid are obtained by Dupin's theorem.  The ellipsoid belongs to a triple orthogonal family of surfaces, formed by the ellipsoid and two hyperboloids (one of one leaf and the other of two leaves).

 For more recent and historical developments of principal lines on surfaces see  \cite{soto1993}, \cite{soto2007}, \cite{garcia2009},  \cite{soga2016} and \cite{soto2021}.
  This work is organized as follows. In section \ref{sec:2}
	we recall the basic properties of the Minkowski 3-space and principal lines.
 In section \ref{sec:3} we describe the global behavior of principal lines in the ellipsoid.
 In section \ref{sec:4} we will describe the topological equivalence of the principal configuration of the ellipsoid.
 In section \ref{sec:5} we will show that the geometric inversion in Minkowski 3-space preserves lines of curvature.
  In section \ref{sec:6} we obtain a triple orthogonal system of quadrics. Finally, in
    section \ref{sec:7} the focal set of the ellipsoid is analyzed.
\section{Preliminaries}\label{sec:2}


The Minkowski 3-space $\mathbb{R}^{2,1}=(\mathbb{R}^3,\langle,\rangle)$ is the vector   space $\mathbb{R}^3$ endowed with the inner product  $\langle u,v\rangle=u_1 v_1+u_2v_2-u_3v_3$, where  $u=(u_1,u_2,u_3)$ and $v=(v_1,v_2,v_3)$. The norm is $\Arrowvert v\Arrowvert=\sqrt{\arrowvert \langle v,v\rangle\arrowvert}$.\\

The vector product $u\times v$, is a vector such that $\langle u \times v,u\rangle=\langle u\times v,v\rangle=0$.  Then
$$u\times v=\left|\begin{matrix}
i & j & -k\\
u_1 & u_2 & u_3\\
v_1 & v_2 & v_3
\end{matrix}\right|.$$ 

\noindent A vector $v$ is said to be

\begin{itemize}
	\item  spacelike,  if $\langle v,v\rangle >0$ or $v=0$,
	\item timelike, if  $\langle v,v\rangle <0$,
	\item lightlike, if $\langle v,v\rangle =0$ and $v\neq 0$.
\end{itemize} 

\noindent  A plane is called spacelike (resp. timelike, lightlike), if the normal vector is  timelike (resp. spacelike,  lightlike). \\

\noindent A regular curve is spacelike (resp. timelike,  lightlike) if the tangent vector is spacelike (resp. timelike,  lightlike). A smooth surface is called spacelike (resp. timelike) if the tangent planes are spacelike (resp. timelike).\\

\noindent Let $\alpha:M\rightarrow \mathbb{R}^{2,1}$ be a $C^r$ ($r\geq 4$) immersion of a smooth and oriented surface $M$ of dimension two in $\mathbb{R}^{2,1}$. Let $X(u,v):\mathbb{R}^2\rightarrow M$ be a local parametrization. The first fundamental form is
\[	I=E du^2+2F du dv+G dv^2, \]
where $E=\langle Xu,Xu\rangle$,  $F=\langle Xu,Xv\rangle$ and $G=\langle Xv,Xv\rangle$. \\

\noindent Given  $p \in M$,   if  $det(I_p)=E G-F^2$  is positive (resp. negative), the surface is spacelike or Riemannian (resp. timelike or Lorentzian) in the point $p$. This is equivalent to say  that tangent plane is spacelike or timelike. The metric induced on $M$ can be degenerate; this happens at the points $p$ on $M$ where the tangent space $TM_p$ is lightlike, or equivalently that $det(I_p)=E G-F^2=0$. We call this set of points the \textit{tropic} and will be denoted by $LD$ (Locus of Degeneracy). \\

 \noindent On a spacelike (resp. timelike) surface, we define the \textit{  Gauss map} $$N(u,v)=\displaystyle \epsilon \cdot\frac{\alpha_u \times \alpha_v}{\Arrowvert \alpha_u \times \alpha_v\Arrowvert}(u,v)$$ 
 such that $N:M \rightarrow H^{2,1}$ with $\epsilon=1$ (resp. $N:M \rightarrow S^{2,1}$ with $\epsilon=-1$), where $H^{2,1}=\{(x,y,z)\in \mathbb{R}^3:x^2+y^2-z^2=-1\}$ and $S^{2,1}=\{(x,y,z)\in \mathbb{R}^3:x^2+y^2-z^2=1\}$.\\

\noindent  The sign $\epsilon=\pm 1$ is only necessary to define the base positively oriented  $\{\alpha_u,\alpha_v,N\}$ in all over the surface (except in the tropic), this is,  $det(\alpha_u,\alpha_v,N)=\displaystyle\frac{\epsilon}{\Arrowvert \alpha_u \times \alpha_v\Arrowvert}\langle \alpha_u\times \alpha_v,\alpha_u \times \alpha_v\rangle>0$,\cite[page 50]{rafael2014differential}.  \\
 
 \noindent The second fundamental form is \[II=edu^2+2f du dv+g dv^2,\]
 where $e=\langle X_{uu},N\rangle$,  $f=\langle X_{uv},N\rangle$ and  $g=\langle X_{vv},N\rangle$.\\

\noindent  The mean curvature $H$ and Gauss curvature $K$ are defined by $$H=\frac{Eg+Ge-2Ff}{2(EG-F^2)}\,\,\,and\,\,\,K=\frac{eg-f^2}{EG-F^2},$$
and the principal curvatures $k_1$ and $k_2$ are defined by $$k_1=H+\sqrt{H^2-K}\,\,\,and\,\,\,k_2=H-\sqrt{H^2-K}.$$

\noindent In general,  a surface  $M\subset  \mathbb{R}^{2,1}$   has a Riemannian part and a Lorentzian part. On the Riemannian part, $dN_p$ does have real eigenvalues; on Lorentzian part, $dN_p$ does not always have real eigenvalues. These eigenvalues are the principal curvatures $k_1$ and $k_2$ in each point and the respective eigendirections of  $dN_p$ are called principal directions and they define two line fields $\mathcal{L}_1$ and $\mathcal{L}_2$ mutually orthogonal in $M$. They  are determined by non-zero vectors 
on $T_p(M)$ which satisfy the implicit differential equation
\begin{equation}\label{eq:principalG}
\left(Fg-Gf\right)dv^2+\left(Eg-Ge\right)du dv+\left(Ef-Fe\right)du^2=0.
\end{equation}
 
The integral curves of the equation \eqref{eq:principalG} are called \textit{ lines of curvature  or principal lines}. The families of principal lines $\mathcal{F}_1$ and  $\mathcal{F}_2$ associated with $\mathcal{L}_1$ and $\mathcal{L}_2$, respectively,  are called principal foliations of $M$. An umbilic point is defined as a point where $II=c I$ for some constant $c$. It is    called a spacelike (resp. timelike) umbilic point  when it is    on Riemannian (resp. Lorentzian) part of $M$. The set of umbilic points is denoted by $\mathcal{U}$.\\

The map $N$ is not defined on the tropic, but since the equation (\ref{eq:principalG}) is homogeneous,  we can multiply the coefficients of (\ref{eq:principalG}) by $\Arrowvert \alpha_u \times \alpha_v\Arrowvert$. Let $L_1=\Arrowvert \alpha_u \times \alpha_v\Arrowvert \left(Fg-Gf\right)$, $M_1=\Arrowvert \alpha_u \times \alpha_v\Arrowvert \left(Eg-Ge\right)$ and $N_1=\Arrowvert \alpha_u \times \alpha_v\Arrowvert \left(Ef-Fe\right)$.  So, \textit{the equation of curvature lines (or principal lines)} can be extended to the tropic by
\begin{equation}\label{eq:principalG2}
L_1 dv^2+M_1 dvdu+N_1du^2=0.
\end{equation}

\noindent   The tropic $LD=(EG-F)^{-1}(0)$ is generically a curve that is solution of the equation \eqref{eq:principalG2}, \cite[Lemma 1.31]{tejada2018}.    The discriminant  of the equation \eqref{eq:principalG2}, define the set of points where it determines a unique direction or an umbilic point, the first set is denoted by $LPL$ (Ligthlike Principal Locus). On the Riemannian part $LPL=\emptyset$, and on the Lorentzian part the set $LPL$ is generically a curve that divide locally the surface  into two regions, in one of them there are no  real principal directions and in the other there are two real principal  directions at each point, \cite{izumiyatari2010}.\\
 
\begin{defi}
    The quintuple $\mathbb{P}_M=\{\mathcal{F}_1,\mathcal{F}_2,\mathcal{U},LD,LPL\}$  is called  the \textbf{principal configuration} of $M$, or rather of the immersion $\alpha$ of $M$ in $\mathbb{R}^{2,1}.$
\end{defi}
\begin{defi}
    Two principal configurations $\mathbb{P}_{M_1}$ and  $\mathbb{P}_{M_2}$ are   $C^0$-principally equivalent if there exists a homeomorphism $h:M_1\to M_2$
    which is a topological equivalence between them, i.e., h sends principal foliations, umbilic set, LD and LPL of $M_1$ in the correspondent of $M_2$.
\end{defi}

\begin{remark}\label{umbilicpoints}
The umbilic points can also be seen as the points where $L_1=M_1=N_1=0$.
\end{remark}
\begin{remark}\label{principalline}
	A smooth curve $c$ is a \textit{principal line}, if  this curve satisfies the equation (\ref{eq:principalG2}) and there are no umbilic points on $c$. 
\end{remark}  

\begin{remark}\label{pricipalchart}
   Let  $X(u,v)$ be a local parametrization of $M$. If  $F=f=0$ then $L_1=N_1=0$ and    $(u,v)$  is a principal curvature coordinate system. It is  called a \textit{ principal chart.} 
\end{remark}  


\subsection*{Triply orthogonal system (see \cite{klingenberg, Struik}).}

In this subsection, it will be introduced a triple orthogonal systems of surfaces in the Minkowski space $\mathbb{R}^{2,1}$.
\begin{defi}
	A triply orthogonal system of surfaces is a differentiable map $X:W\longrightarrow \mathbb{R}^{2,1}$, defined on an open set $W \subset \mathbb{R}^{2,1}$, satisfying:
	\begin{itemize}
		\item[a)] The linear map $dX_{(u,v,w)}:T_{(u,v,w)}\mathbb{R}^{2,1}\longrightarrow T_{X(u,v,w)}\mathbb{R}^{2,1}$ is bijective for all $(u,v,w)\in W$.
		\item[b)] $\langle X_u,X_v\rangle=\langle X_u,X_w\rangle=\langle X_w,X_v\rangle=0$.
	\end{itemize}
\end{defi}

 Let $p=(u_0,v_0,w_0)\in W$.  Consider the three surfaces
$$(u,v)\longmapsto X(u,v,w_0)$$
$$(u,w)\longmapsto X(u,v_0,w)$$
$$(v,w)\longmapsto X(u_0,v,w),$$
  we denote these surfaces by $M_{w_0}$, $M_{v_0}$ and $M_{u_0}$, respectively. They are regular surfaces by  the condition a).\\

 Notice that by condition b), $F=0$ on each of them. Furthermore, $X_{w}(u,v,w_0)$ is normal to $M_{w_0}$ at $(u,v,w_0)$ (similarly to other two surfaces) and differentiating,
 $$\langle X_u,X_v\rangle_w=\langle X_u,X_w\rangle_v=\langle X_w,X_v\rangle_u=0.$$

\noindent Therefore, $$\langle X_{uv},X_w\rangle=\langle X_{uw},X_v\rangle=\langle X_{vw},X_u\rangle=0,$$
which means that $f=0$ on each of the surfaces. By remark (\ref{pricipalchart}), we may conclude that:
 
\begin{thm}\label{Dupin}
The coordinate curves on a surface belonging to  a triply
orthogonal system in a Minkowski tridimensional space are   principal curvature lines.
\end{thm}



\section{The Ellipsoid in the Minkowski space}\label{sec:3}

Consider the family of surfaces $$\displaystyle \mathbb{F}_{u}=\{(x,y,z): \frac{ x^2}{a^2- u}+\frac{ y^2}{b^2- u}+\frac{ z^2}{c^2+ u}=1\}$$  $$\displaystyle \mathbb{G}_v=\{(x,y,z): \frac{ x^2}{a^2-v}+\frac{y^2}{b^2-v}+\frac{z^2}{c^2+v}=1\}$$ 

$$\displaystyle \mathbb{H}_w=\{(x,y,z): \frac{x^2}{a^2-w}+\frac{y^2}{b^2-w}+\frac{z^2}{c^2+w}=1\}$$ 
\noindent where $a>b>0$ (the case $b>a>0 $ is similar) and $c>0$.\\

Let $U_E:=\{(u,v,w) \in  (-c^2,b^2)\times (b^2,a^2) \times (-c^2,b^2)\}$. For $(u,v,w) \in U_E$,  $F_u$, $H_w$  are ellipsoids and  $G_v$ is a  hyperboloid of one leaf.\\

\begin{thm}\label{ortogonalE}
The surfaces $\mathbb{F}_u$, $\mathbb{G}_v$ and $\mathbb{H}_w$ define a triple orthogonal system for $(u,v,w)\in U_E$, $u\neq w$. 
\end{thm}
	
\begin{proof}
Solving the system  below in the variables $\{x ,y ,z \}$
  \begin{equation*}
  \begin{split}
  &\frac{x^2}{a^2-u}+\frac{y^2}{b^2-u}+\frac{z^2}{c^2+u}-1=0\\
  &\frac{x^2}{a^2-v}+\frac{y^2}{b^2-v}+\frac{z^2}{c^2+v}-1=0\\
  &\frac{x^2}{a^2-w}+\frac{y^2}{b^2-w}+\frac{z^2}{c^2+w}-1=0,	
  	\end{split}
  	\end{equation*}
it is obtained in the positive octant:
    \begin{equation*}
    \begin{split}
   x(u,v,w)= &  \sqrt{\frac{(a^2-u)(a^2-v)(a^2-w)}{(a^2-b^2)(a^2+c^2)}}\\
   y(u,v,w)= &   \sqrt{\frac{-(b^2-u)(b^2-v)(b^2-w)}{(a^2-b^2)(b^2+c^2)}}\\
   z(u,v,w)= &  \sqrt{\frac{(c^2+u)(c^2+v)(c^2+w)}{(a^2+c^2)(b^2+c^2)}}	\cdot
    \end{split}
   \end{equation*}
 A long and straightforward calculations show that   \begin{equation}\label{eq:chart1E}
   X(u,v,w)=(x(u,v,w),y(u,v,w),z(u,v,w))
   \end{equation} satisfies   $\langle X_u,X_v\rangle=\langle X_u,X_w\rangle=\langle X_v,X_w\rangle=0$. Moreover, 
   {\small
   \[\textrm{det}(DX(u,v,w))=\frac{(u-v)(u-w)(v-w)}{8 x(u,v,w) y(u,v,w) z(u,v,w) (a^2-b^2)(a^2+c^2)(b^2+c^2)}\neq 0.\]
   }
\end{proof}

Since $\{F_u, G_v, H_w\}$  is a triple   orthogonal system, these surfaces intersect along their curvature lines. The curvature lines   can be obtained  globally by symmetry in relation to the coordinates planes.\\

\begin{figure}[H]
	\begin{center}
		\def\svgwidth{0.7\linewidth}
 	\includegraphics[width=0.7\textwidth,angle=0]{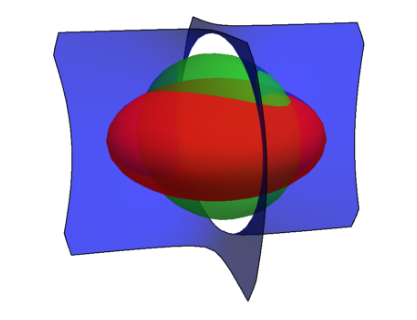}
		\caption{Triply orthogonal system defined by two   ellipsoids and one hyperboloid of one leaf.}
	\label{stoelipsoide}
	\end{center}
\end{figure}

 Now, we fixed $w$ and defined the ellipsoid $\mathbb{E}_w=\mathbb{H}_w$ with $(u,v,w) \in U_E$, so we have that the principal lines on $\mathbb{E}_w$ are the intersection curves, with the hyperboloid of one leaf $\mathbb{F}_u$  and with the other ellipsoid $\mathbb{G}_v$ (See Fig. \ref{stoelipsoide}).\\

In each octant, we have that for $E_w$ the parametrization (\ref{eq:chart1E}) is a principal chart, i.e., $f=F=0$. So, the principal lines are $u=constant$ and $v=constant$, and these curves are exactly the intersection between surfaces. \\

\begin{remark}\label{stoE}
  The triply orthogonal system of quadratic surfaces  in the Euclidean space is make up by an ellipsoid, a hyperboloid of one leaf and a hyperboloid of two leaves \cite[Chapter 2]{garcia2009}.
See also \cite[Chapter 7]{stachel} for more details about the geometric properties of confocal quadrics.
\end{remark}

To complete the description, the principal configuration is necessary to analyze the curves of the intersections of the ellipsoid with the coordinates planes. Without loss of generality, we do $w=0$, i.e., we analyze 
\[  \mathbb{E}_0=\{(x,y,z): \frac{x^2}{a^2}+\frac{y^2}{b^2}+\frac{z^2}{c^2}=1\}\]
with $a>b>0$ and $c>0$ (it is allowed $a=c$ or $b=c$).\\

The parametrization below is inspired in the Euclidean case. See also Section \ref{sec:6} where a global parametrization will be obtained in a triple orthogonal system of quadrics.
\begin{lema}\label{globalchart}
The parametrization 
	\begin{equation}\label{chart2E}
	\aligned
	X(u,v)=&\left(a\cos \left( u \right) A(v) ,b \sin(u) \sin(v),c B(u) \cos(v) \right)\\
	A(v)=&\sqrt{A_1 \cos^2(v)+\sin^2(v)}, B(u)=\sqrt {B_1 \cos^2(u)+\sin^2(u)}
	\endaligned
	\end{equation}		with $(u,v)\in U_1=[0,\pi]\times [0,2\pi]$  or $(u,v)\in U_2=[0,2\pi]\times [0,\pi]$,	where $A_1=\frac{a^2-b^2}{a^2+c^2}$ and $B_1=\frac{b^2+c^2}{a^2+c^2}$,
	defines a principal chart $(u,v)$ on the ellipsoid $\mathbb{E}_0$.
\end{lema}
 
\begin{proof} Calculating of the coefficients of the first and second fundamental form, we have

\begin{align*}
    E=&-{\frac { \left(  \left( {a}^{2}-{b}^{2} \right)    \cos^2   u
  -{a}^{2} \right)  \left(  \left( {a}^{2}-{b}^{2}\right)    \cos^{2}u+ \left( {b}^{2}+{c}^{2} \right)   \cos^{2}v-{a}^{2}-{c}^{2} \right) }{ \left( {a}^{2}-{b}^{2} \right)    \cos^{2}u-{a}^{2}-{c}^{2}}}\\
    F=&0\\
    G=&-{\frac { \left(  \left( {a}^{2}-{b}^{2} \right)  \cos^2 u+ \left( {b}^{2}+{c}^{2} \right)  \cos^2 v-{a}^{2}-{c}^{2} \right)  \left( 
 \left( {b}^{2}+{c}^{2} \right)  \cos^2 v-{c}^{2} \right) }{ -\left( {b}^{2}+{c}^{2} \right) 
\cos^2 v+{a}^{2}+{c}^{2}}}
\\
    e=&\frac{abc \left(  \left( {a}^{2}-{b}^{2} \right)  \cos^2 u+ \left( {b}^{2}+{c}^{2} \right) \cos^2 v-{a}^{2}-{c}^{2} \right) ^{2}  }{\sqrt { \left( {b}^{2}+{c}^{2} \right)  \cos^2 v-{c}^{2}-{a}^{2}}    \left(  \left( {a}^{2}-{b}^{
2} \right)  \cos^2 u-{a}^{2}-{c}^{2}
 \right) ^{{\frac{3}{2}}}}\\
    f=&0\\
    g=&\frac{b \left(  \left( {a}^{2}-{b}^{2} \right) \cos^2 u+ \left( {b}^{2}+{c}^{2} \right) \cos^2 v-{a}^{2}-{c}^{2} \right) ^{2} a c}{\left(  \left( {b}^{2}+{c}^{2} \right)  \cos^2 v-{c}^{2}-{a}^{2}\right) ^{\frac 32}    \sqrt { \left( {a}^{2}-{b}^{
2} \right) \cos^2 u-{a}^{2}-{c}^{2}}}.
\end{align*}

 Since that, $F=f=0$ then $L_1=N_1=0$ in the equation (\ref{eq:principalG2}).  So, by Remark \ref{pricipalchart}, we have that $X$ defines a parametrization by principal lines, i.e.,  $(u,v)$ is a principal chart. \\
 
 The parametrization $(X,U_1)$ (resp. $(X,U_2)$) cover all the ellipsoid and is smooth, except in the curves $X(0,v)$ and $X(\pi,v)$ (resp. $X(u,0)$ and $X(u,\pi)$).
\end{proof}
 
\begin{pro}\label{ProE}
	On the ellipsoid $\mathbb{E}_0$, we have that:
	\begin{itemize}
		\item[a.] The curves $c_x=\{(x,y,z):x=0\}\cap \mathbb{E}_0$ and $c_z=\{(x,y,z):z=0\}\cap \mathbb{E}_0$ are principal lines. 
		\item[b.] $\mathbb{E}_0$ has exactly four spacelike umbilic points,  $$\left(\pm a\sqrt{\frac{a^2-b^2}{a^2+c^2}},0,\pm c\sqrt{\frac{b^2+c^2}{a^2+c^2}}\right).$$
		\item[c.] The umbilic points are of type $D_1$.
		\item[d.] The curve $c_y=\{(x,y,z):y=0\}\cap \mathbb{E}_0$ is the union of principal lines. Moreover, these are the separatrices of the umbilic points.
		\item[e.] The tropic is composed by two disjoint regular closed curves. Moreover, these curves are principal lines.
		\item[f.] The principal lines are globally defined, i.e., the set $LPL=\emptyset$.
	\end{itemize}
\end{pro}

 

   	\begin{figure}[H]
 
 	\centering	
 \includegraphics[width=0.7\textwidth,angle=0]{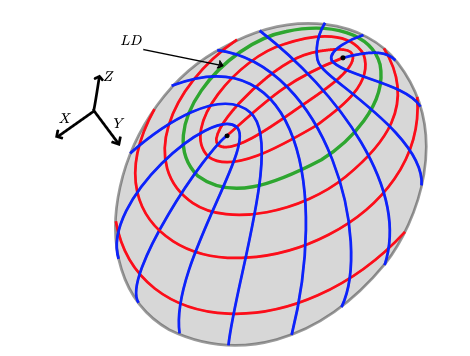}
  \caption{	Principal lines on the Ellipsoid in the Minkowski space. Parameters $a=2.0,\,\, b=1.5,\,\, c=2.2$.  \label{fig:ellipsoid0} }
 	\end{figure}
 
\begin{proof} 
\noindent \textbf{a)\;} Consider the principal chart $(X,U_1)$ (resp. $(X,U_2)$) given by Lemma \ref{globalchart}. We have $c_x=X(\frac \pi 2,v)$ (resp. $c_z=X(u,\frac \pi 2)$).\\

The principal chart $(X,U_1)$ (resp. $(X,U_2)$) is smooth, except in the curves $X(0,v)$ and $X(\pi,v)$ (resp. $X(u,0)$ and $X(u,\pi)$), but this curves not intersect with $c_x$ (resp. $c_z$). Therefore, $c_x$ (resp. $c_z$) is a principal line of the ellipsoid.\\

\noindent \textbf{b)\;} Consider the parametrization, 
    \begin{equation} \label{chart3E} X(u,v)=\left(a u, b v, \pm c \sqrt{1-u^2-v^2} \right).\end{equation}
	Then the differential equation of principal lines (\ref{eq:principalG2}) with $X$ is   
 
	\begin{equation}\label{eqchart3E}\aligned
	E(u,v,du:dv)&= -uv(a^2+c^2)du^2+ uv(b^2+c^2)dv^2 \\
 &+(u^2(a^2+c^2)-v^2(b^2+c^2)-a^2+b^2)dudv 
 =0.
 \endaligned
	\end{equation} 
		
	We have that
		$L_1=N_1=0$ when $u=0$ or $v=0$. If $u=0$ then $M_1=-v^2(b^2+c^2)-a^2+b^2\neq 0$. If $v=0$, we have that $M_1=u^2(a^2+c^2)-a^2+b^2=0$ if and only if $$u_0=\pm\sqrt{\frac{a^2-b^2}{a^2+c^2}}.$$
		
	So,   there are four umbilic points. 	Moreover, $$(E G-F^2)(u_0,0)=\frac{b^4(a^2+c^2)}{b^2+c^2}>0,$$ and then the umbilic points are in the Riemannian part of $E_0$, i.e, they are spacelike umbilic points.\\

\noindent \textbf{c)\;} For completeness, it will be included a detailed sketch of proof.
We take $ p=\frac{dv}{du}$ in the equation \eqref{eqchart3E}, so 	$$F(u,v,p)=-uv(a^2+c^2)+(u^2(a^2+c^2)-v^2(b^2+c^2)-a^2+b^2)p+uv(b^2+c^2)p^2=0.$$
  Under the hypothesis, the implicit surface $F^{-1}(0) $ is a regular surface, contain the projective line, and is topologically a cylinder. The map $\pi:F^{-1}(0)\to \mathbb{R}^2$, $\pi(x,y,p)=(x,y)$ is a ramified double covering and $\pi^{-1}(u_0,0)$
    is the projective line parametrized by $[du:dv]$.
\noindent The umbilic point $ P_1=\left( a\sqrt{\frac{a^2-b^2}{a^2+c^2}},0, c\sqrt{\frac{b^2+c^2}{a^2+c^2}}\right)$ has coordinates $u_0=\sqrt{\frac{a^2-b^2}{a^2+c^2}}$ and $v=0$. 
\noindent The Lie-Cartan line field associated to the implicit equation $F(u,v,p)=0$ is $Y=(F_p,pF_p,-(F_x+pF_y))$ on the surface  $M=F^{-1}(0)$, $M\subset \mathbb{R}^2\times \mathbb{RP}^1$. The solutions of the implicit differential equation $F(u,v,p)=0$ are the projections of the integral curves of $Y$. See \cite{BF} and \cite{garcia2009}.

We have that 
				
	$$Y(u_0,0,p)=\left(0,0,-\sqrt{\frac{a^2-b^2}{a^2+c^2}}p(b^2p^2+c^2p^2+a^2+c^2)\right)=(0,0,0)$$ if and only if $p=0$. \\
				
	Moreover, the eigenvalues of  $DY(u_0,0,0)$ are, $$\lambda_1=2\sqrt{\frac{a^2-b^2}{a^2+c^2}}(a^2+c^2),\;\;\;\lambda_2=-2\sqrt{\frac{a^2-b^2}{a^2+c^2}}(a^2+c^2).$$
		Therefore, $(u_0,0,0)$	is a hyperbolic saddle point of $Y$. To complete the analysis, it is also necessary to consider the chart $q=du/dv$ in the equation \eqref{eqchart3E} to obtain an implicit surface $G(u,v,q)=0$. Now the Lie-Cartan vector field
  is $Z=(qG_q,G_q,-(qG_u+G_v))$. We have that $Z(u_0,0,q)\ne 0$. Gluing the phase portraits of $Y$ and $Z$ near the projective line  $[du:dv]$ we obtain a line field on the cylinder 
  with a unique hyperbolic singular point. The projections of the leaves (integral curves of $X$ and $Y$) are the principal lines of the ellipsoid  near the umbilic point. 
  See Fig. \ref{fig:cilindrosela}.

	\begin{figure}[H]
	\begin{center}
\includegraphics[width=0.9\textwidth,angle=0]{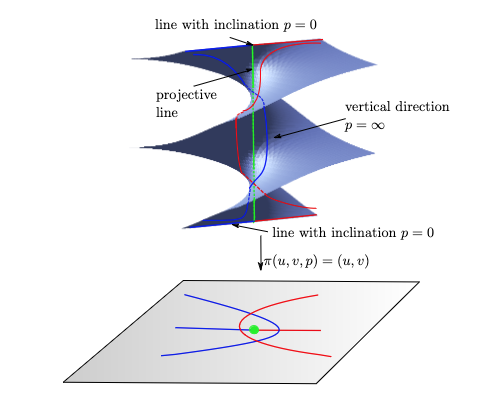}
		\caption {Implicit surface $F(u,v,p)=0$ (cylinder) and a ra\-mi\-fied double covering $\pi$ with $\pi^{-1}(u_0,0)$ being the projective line. The top and bottom lines with inclination $p=0$ are identified.  \label{fig:cilindrosela}}
\end{center}
\end{figure}
  Therefore, the umbilic point $P_1$ is Darbouxian of type $D_1$ (see also \cite{garcia2009}). By symmetry,  all the other umbilic points are also of type $D_1$. \\
\noindent \textbf{d)\;} Using the parametrization (\ref{chart3E}), a curve $c_y$ satisfies the equation the principal lines (\ref{eqchart3E}). Furthermore, the umbilic points are on $c_y$, so this curve  is a union of principal lines and the umbilic points. Since the umbilic points are $D_1$, we obtain the result as stated.\\
	
\noindent \textbf{e)\;} Using the chart defined by equation (\ref{chart2E}) with $(u,v)\in U_1$,
	we have that
	\begin{equation*}
	\begin{split}
	E G-F^2=&\left((b^2+c^2)\cos^2(v)-c^2\right)\left((a^2-b^2)\cos^2(u)-a^2\right)\\
         	&\left((a^2-b^2)\cos^2(u)+(b^2+c^2)\cos^2(v)-a^2-c^2\right)^2=0.
	\end{split}
	\end{equation*} 
	 if, and only if, 
	 $ v_1=\arccos\left(\frac{c}{\sqrt{c^2+b^2}}\right)$ or $v_2=\arccos\left(-\frac{c}{\sqrt{c^2+b^2}}\right)=\pi-v_1$.\\
	 
	The tropic is the union of the closed curves $c_1(u)=X(u,v_1)$ and $c_2(u)=X(u,v_2)$. As $v=constant$ and $(X,U_1)$ is a principal chart, then $c_1$ and $c_2$ are principal lines.\\
	
\noindent \textbf{f)\;} Since    the parametrization (\ref{chart2E}) is defined globally on $\mathbb{E}_0$  and defines a principal chart, it follows that   $L_1=N_1=0$ and $LPL=M_1^2\geq 0$. Therefore,  the principal lines are globally defined. 
    \end{proof}

\noindent From Proposition \ref{ProE} and Theorem \ref{ortogonalE}, the principal configuration on the Ellipsoid $\mathbb{E}_0$ is as shown in  Fig. \ref{fig:ellipsoid0}. All principal lines are closed, except of four open arcs  that are the connections between umbilic points.  They are called  the \textit{umbilic separatrices}.
 
\subsection*{Confocal and orthogonal family of conics}
   Performing the change of   coordinates by $u=\sqrt{\frac{b^2+c^2}{a^2+c^2}}\; x$ and $v=y$, then equation (\ref{eqchart3E})  is given by 
   \[-xydx^2+\left(x^2-y^2-\lambda^2\right)dx dy+xydy^2=0\]
    with $\lambda^2=\frac{a^2-b^2}{b^2+c^2}$. The coordinates axes,  the family of ellipses
    \[u(t)=R \cos(t), v(t)=r \sin(t), R^2=r^2+\lambda^2\]
     and the family of  hyperbolas 
 \[u(t)=R \cosh(t), v(t)=r \sinh(t), R^2-r^2=\lambda^2\]
   are the solutions of the differential equation above. This  is similar to  Euclidean case, see \cite{garcia2009}.
   
       \begin{figure}[H]
   	
   	\centering
   	\includegraphics[width=0.9\textwidth,angle=0]{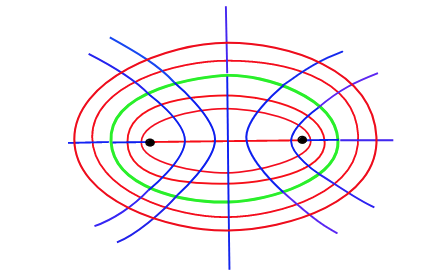}	
   	\caption{	\label{elliphiper} Confocal and orthogonal family of conics.The tropic is shown in green and is parametrized by $\cos^2 v=c/(b+c)$.}
   \end{figure}

\subsection*{\bf Horizontal ellipsoid of revolution.} When $a=c$ or $b=c$, we have four spacelike umbilic points of type $D_1$, while with the  Euclidean scalar product only have two umbilic points of center type.

\subsection* {\bf Vertical ellipsoid of revolution and Euclidean sphere.}
  When $a=b$ and $c>0$, the parametrization (\ref{chart2E}) is reduced to
  \[X(u,v)=(a\sin(v)\cos(u),a\sin(v) \sin(u),c \cos(v)).\]
 \noindent The equation of principal lines is
 \[a^2c (a^2 \cos^2(v)+c^2\cos^2(v)-a^2-c^2)^4   du dv=0.\]
Therefore,  the principal lines are $u=constant$ and $v=constant$. We have only two spacelike umbilic points $(0,0,\pm c)$ of type center. \\
  
 \noindent On ellipsoid of revolution with $a=b$ and $c\neq a$, the principal lines are the same in the two geometries (Euclidean and Lorentzian).\\
 
 \noindent With the Euclidean scalar product the Eucidean sphere is a umbilic surface, while with the  Lorentzian scalar product the Euclidean sphere has only two spacelike umbilic points of type center.

\subsection*{\bf Umbilic surfaces.} The umbilic surfaces with Euclidean inner scalar are the Euclidean sphere and planes, while with Lorentzian inner scalar the umbilic surfaces are planes,  the vertical hyperboloid of one leaf $\mathbb{S}^2_1(c,r) =\{p\in\mathbb{R}^{2,1}:\langle p-p_0,p-p_0\rangle=r^2\}$  and vertical hyperboloid of two leaves   $\mathbb{H}^2_1(c,r)=\{p\in\mathbb{R}^{2,1}:\langle p-p_0,p-p_0\rangle=-r^2\}$, see   \cite[page 191]{ivo} and \cite[page 85]{rafael2014differential}.
\begin{remark}
    For the study of geodesics on an ellipsoid in the Minkowski space $\mathbb{R}^{2,1}$ see \cite{taba_2007}. The analysis of umbilic points in smooth surfaces in $\mathbb{R}^{2,1}$ of the form $f_{\varepsilon}(x,y,z)=x^2/a^2+y^2/b^2+z^2/c^2+h.o.t=\varepsilon$ was developed in \cite{tari_2017}.
\end{remark}
\section{Topological equivalence of principal foliations}\label{sec:4}
In this section we will obtain that the principal configurations of the ellipsoids of three distinct axes are all principal topologically equivalent. The  Euclidean case was established by J. Sotomayor \cite{soto1993}.

\begin{pro}\label{prop:canonical}
   Consider an ellipsoid $E(x,y,z)= ax^2+by^2+cz^2+2dxy+2exz+2fyz+gx+hy+kz+l=0$.
   Then there exists an isometry   $h:\mathbb{R}^{2,1}\to \mathbb{R}^{2,1}$ such that $E(h(u,v,w))=\lambda_1u^2+\lambda_2v^2+\lambda_3w^2=1$, with $\lambda_i>0$ for $(i=1,2,3).$

\end{pro}
\begin{proof}    
The rotation group  of $\mathbb{R}^{2,1}$ is $SO(2,1)$ of dimension 3 and is generated by the Euclidean and Hyperbolic rotations defined by:  

\begin{equation*}
    \aligned
    R(u,v,w)&=( u\cos\theta  +v\sin\theta  , -u\sin\theta  +v\cos\theta  ,w)\\
    S(u,v,w)&=(u \cosh\alpha  +w\sinh\alpha  ,v,  u\sinh\alpha  +w \cosh\alpha  )\\
    T(u,v,w)&=( u, v\cosh\beta  +w \sinh\beta  ,  v\sinh\beta  +w \cosh\beta)
    \endaligned
\end{equation*}

The quadric form $q(x,y,z)= ax^2+by^2+cz^2+2dxy+2exz+2fyz$ is positive definite when  one of the following conditions holds.
   \begin{align*}
       &a> 0,\; ab-d^2>0,\; abc - af^2 - be^2 - cd^2 + 2de f=\Delta>0,\\
      &b>0,\; bc-f^2>0,\; \Delta>0,\\
      &c>0,\; ac-e^2>0,\; \Delta >0.
   \end{align*} 
    
   In this case the eigenvalue problem
   \[
  {\rm det} \left( \begin{matrix} a-x & d & e\\
   d&b-x & f\\
   e &f&c+x\end{matrix}\right)=0
   \]
   has three real eigenvalues $x_1\leq x_2\leq x_3$ and the correspondent eigenvectors $e_1, e_2, e_3$ are orthonormal relative to the Minkowski inner product. Therefore, one of the eigenvectors, say $e_3$, 
   is timelike and the other two $\{e_1,e_2\}$   are spacelike.
   There exists an isometry $h$ (composition of hyperbolic rotations) such that $h(0,0,-1)=e_3$.
   In the new coordinates we have that
   $q_1(u,v,w)=q(h(u,v,w))=a_1u^2+b_1v^2+c_1w^2+d_1 uv$.

Also, there exists an isometry $H$ (Euclidean rotation) such that 
\[  q_1(H(u_1,v_1,w_1))=a_2u_1^2+b_2 v_1^2+c_1w_1^2,\;\; \;a_2>0,\; b_2>0, \; c_1>0.\]
Finally, with a translation we obtain the result stated.
    \end{proof}

\begin{remark} In general, a hyperboloid is not isometric to one given in a diagonal form. The classification of conics in Minkowski plane  is carried out in  \cite{ruas1987}.
    
\end{remark}
\begin{thm}
\label{th:equivalence}
    Consider the set of ellipsoids $\mathcal{E}$ with three distinct axes in the space of quadrics $\mathcal{Q}$ of $\mathbb{R}^3$. Then the principal configurations of any two elements of $\mathcal{E}$ are principal topologically equivalent.
    
\end{thm}

\begin{proof} The principal configuration of an ellipsoid  with three different axes in the diagonal form $x^2/a^2+y^2/b^2+z^2/c^2=1$ has the following properties.
\begin{itemize}
    \item[i)] There are four umbilic points of Darbouxian type $D_1$.
    \item[ii)] The set LD is the union of two regular curves.
    \item[iii)] The set LPL is empty.
    \item[iv)] The principal foliations $\mathcal{F}_1$ and $\mathcal{F}_2$ have all leaves closed, with the exception of the umbilic separatrices. See Fig. \ref{fig:ellipsoid0}.
    
\end{itemize}

    The construction of the topological equivalence can be done explicitly using the method of canonical regions defined by the union of two topological disks  and a cylinder;  the boundary being the tropics. See \cite{soto1982a} and Fig. \ref{fig:confelip}.
       \begin{figure}[H]
   	\centering
   	\includegraphics[width=1.0\textwidth,angle=0]{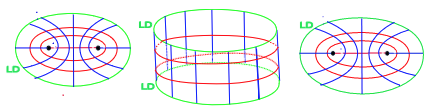}	
   	\caption{	\label{fig:confelip}  Decomposition of the ellipsoid in three canonical regions foliated by principal lines; the boundary of each region is formed by    the tropic lines. }
   \end{figure}
By Proposition \ref{prop:canonical},   any ellipsoid is isometric to an ellipsoid in the  diagonal form.
    This ends the proof.
\end{proof}

\section{Geometric Inversion in Minkowski space }\label{sec:5}
In this section, we will show that the principal lines are the same when we consider the  inversion of the surface with respect to a given  point in the space. Recall that the inversion is defined by:
\[ I_q(p)=\frac{p-q}{\langle p-q, p-q\rangle}.\]
 
\begin{pro}\label{prop:inversion} Consider a regular surface $S$ and a point $q\in\mathbb{R}^{2,1}\setminus S$
Let $S_q=I_q(S)$, where $I_q$ is the  inversion with respect to the point $q$. The principal lines on $S$ are the same that on $S_q$.  
\end{pro}

\begin{proof}
Consider the local parametrization $$X(u,v)=(u,v,h(u,v)).$$ 
Calculating the equation of principal lines of $X$, we have that:
\begin{equation*}\label{inversion}
\begin{aligned}
   (&h_{uv}h^2_v-h_{vv}h_uh_v-h_{uv})dv^2+(h_{uu}h^2_v+h_{vv}-h_{vv}h^2_u-h_{uu})dvdu\\
   +&(h_{uv}+h_uh_vh_{uu}-h_{uv}h^2_u)du^2=0.
\end{aligned}
 \end{equation*}
The local parametrization of the inverted surface in the relation to the point   $q=(q_1,q_2,q_3)$, is given by:
\begin{align*}
 \bar{X}(u,v)=&\frac{1}{\langle X(u,v)-q, X(u,v)-q\rangle}(X(u,v)-q)\\
=&\frac{1}{(u-q_1)^2+(v-q_2)^2-(h(u,v)-q_3)^2}  \left(u-q_1,v-q_2,h(u,v)-q_3\right).  
\end{align*}
 
\noindent Calculating the first fundamental form of  $\bar{X}$ it follows that
\begin{align*}
 E=&-\frac{h^2_u-1}{Q_0^2},\,\,\,
F=-\frac{h_uh_v}{Q_0^2},\,\,\,
 G=-\frac{h^2_v-1}{Q_0^2},
\end{align*}
with $Q_0=\langle X(u,v)-q, X(u,v)-q\rangle$. Similarly, we calculate the coefficients of the second fundamental form:

\begin{align*}
 e=&\frac{1}{Q_0^4}\left[ 2(q_1-u)h_u^3+2((q_2-v)h_v+h-q_3)h_u^2+2(u-q_1)h_u+2(v-q_2)h_v\right. \\ 
 & \left.  \,\,\,\,\,+2q_3-2h +(h^2-2hq_3-q_1^2+2q_1u-q_2^2+2q_2v+q_3^2-u^2-v^2)h_{uu} \right]\\
f=&\frac{1}{Q_0^4}\left[ 2(q_1-u)h_vh_u^2+2((q_2-v)h_v^2+(h-q_3)h_v)h_u+h^2h_{uv}-2h h_{uv}q_3\right. \\
 & \left.  \,\,\,\,\,- h_{uv}q_1^2+2h_{uv}q_1u-h_{uv}q_2^2+2h_{uv}q_2v+h_{uv}q_3^2-h_{uv}u^2-h_{uv}v^2\right]\\
 g=&\frac{1}{Q_0^4}\left[2((q_1-u)h_v^2+u-2q_1)h_u+2(q_2-2v)h_v^3+2(h-q_3)h_v^2+2(v-q_2)h_v \right.\\
 &\left.   \,\,\,\,\,+2q_3-2h+(h^2-2hq_3-q_1^2+2q_1u-q_2^2+2q_2v+q_3^2-u^2-v^2)h_{vv} \right].
\end{align*}

\noindent Then, the coefficients of the differential equation of the principal lines   are given by: 
\begin{align*}
 L=Fg-Gf=&\frac{h_uh_vh_{vv}-h_{uv}h^2_v+h_{uv}}{\langle X-q, X-q\rangle^5},\\
M=Eg-Ge=&\frac{h_uh_{vv}-h_{uu}h^2_v+h_{uu}-h_{uv}}{\langle X-q, X-q\rangle^5},\\
L=Ef-Fe=&\frac{h^2_uh_{uv}-h_uh_{uu}h_v-h_{uv}}{\langle X-q, X-q\rangle^5}.
\end{align*}

\noindent So, the differential equation of the principal lines of $\bar{X}$ is exactly (\ref{inversion}). The analysis, with local parametrization $(u,h(u,v),v)$ or $(h(u,v),u,v)$ are analog.

\noindent Therefore, the principal lines of the surface $S$ and of the inverted surface $S_q$ are related by the inversion $I_q$, i.e., if $\gamma(s)$ is a principal line of $S$, then $I_q(\gamma(s))$ is a principal line of $S_q$.
\end{proof}

\section{Triple Orthogonal System in Minkowski space}\label{sec:6}

In this section a global parametrization of a triple orthogonal system of quadrics in the Minkowski 3-space will be established.

Let $Z(u,v,w)= (A(u,v,w),B(u,v,w),C(u,v,w))$ defined by:
{\small 
\begin{equation}\label{eq:sto}
\begin{aligned}
A(u,v,w)=&\cos u\, \cosh   w \sqrt {(\varepsilon\,{n}^{2}+m^2)   \cos^2v   +m^2  \sin^2 v }\; 
  \\
 B(u,v,w)=&  m\sin u \sin v\sinh   w
 \\
 C(u,v,w)=&\cos \left( v \right) \sqrt {{\frac { \left( \varepsilon\,{n
}^{2}  \cos  ^{2}u- m^2 \sin^2  u\right)
   \left(  \varepsilon\,{n}^{2}\cosh ^{2}w+  m^2\sinh^2   w
  \right)   }{\varepsilon {n}^2+{m}^2}}}
\end{aligned}
\end{equation}
}
Here $\epsilon =\pm 1.$
\begin{thm}\label{th:sto}
The map $Z$ defined by equation \eqref{eq:sto} is a triple orthogonal system of quadrics in $\mathbb{R}^{2,1}$ (Minkowski 3-space).
More precisely,  the quadrics are given by:

\begin{equation*}
\label{eq:stoq}
\begin{aligned}
 \mathcal{E}_1: &\; {\frac {{x}^{2}}{  {m}^{2}\cosh^{2}w}
}+{\frac {{y}^{2}}{    {m}^{2}\sinh^2w 
}}+{\frac {{z}^{2} \left( {m}^{2}+\varepsilon{n}^{2} \right) }{    \left(\varepsilon 
{n}^{2}\cosh ^{2}w -   {m}^{2} \sinh ^{2}w \right) {m}^{2}}}=1\\
 \mathcal{E}_2:\; & \frac {{x}^{2} \left( {m}^{2}+\varepsilon{n}^{2} \right) }{{m}^{2} \left( {m}^{2}+ \varepsilon{n}^{
2}  \cos ^{2}v \right) }+{
\frac { \left( {m}^{2}+\varepsilon{n}^{2} \right) {y}^{2}}{\varepsilon {m}^{2}{n}^{2}  
\sin^{2}v }+  \frac {{z}^{2} \left( {m}^{2}+\varepsilon{
n}^{2} \right) }{   \varepsilon{m}^{2}{n}
^{2}\cos^{2}v}}=1
\\
 \mathcal{H}_1: &\; {\frac {{x}^{2}}{ {m}^{2} \cos^{2}u  }}-
{\frac {{y}^{2}}{  {m}^{2} \sin ^{2}u }}+
{\frac {{z}^{2} \left( {m}^{2}+\varepsilon{n}^{2} \right) }{{m}^{2} \left( 
 {m}^{2} \sin ^{2} u + \varepsilon{n}^{2}  \cos^{2}u  \right) }}=1
\end{aligned}
\end{equation*}

\end{thm}

\begin{proof} The map $Z$ defined by equation \eqref{eq:sto} was inspired in \cite{tabanov} where a similar map was obtained in the Euclidean case.
The main idea is to try a parametrization with separation of variables as 
\[ Z(u,v,w)= \left(h_1\cos u\cos w\, a(v),h_2 \sin u\sin v\sinh w,h_3 \,c(u)\,d(w) \,\cos v \right). \]
A long,  and straightforward calculation, using the equation \eqref{eq:sto}, leads to
\[ \langle Z_u,Z_v\rangle =\langle Z_u,Z_w\rangle=\langle Z_v,Z_w\rangle=0. \]
The quadrics defined by equation \eqref{eq:stoq}   was obtained by the method of elimination of variables from the equations
\[ A(u,v,w)-x=0, \;\;B(u,v,w)-y=0, \; \;C(u,v,w)-z=0.\]
\end{proof}

\begin{remark}\label{rem:ellipsoid}
    For $\varepsilon=-1$, $m=\sqrt{a^2 - b^2}$,\, $n=\sqrt{ (b^2 + c^2)(a^2 - b^2)}/\sqrt{a^2 + c^2}$
    and $\cosh w=a/\sqrt{a^2 - b^2}$ it follows that
 $Z_w(u,v)=Z(u,v,w)$ is a parametrization of the ellipsoid $x^2/a^2+y^2/b^2+z^2/c^2=1$.
\end{remark}

\section{Focal set of a surface in the Minkowski space}\label{sec:7}
The focal set of a surface $M$ can be defined as the singular set of the  congruence of lines  given by 
\[ L(u,v,t)=\alpha(u,v)+t N(u,v)\]
where $\alpha$ is a parametrization of $M$ and $N$ is the normal vector to the surface.
Also, the focal set can be seen as the locus of the centers of curvature of the given surface.

\[ \mathcal{F}_i:\; \alpha(u,v)+\frac{1}{k_i(u,v)} N(u,v), \;\;\; (i=1,2).\]
See \cite{cayley1873} and \cite{tari2012}. 

\begin{pro}
\label{prop:focal12}
The focal set $\mathcal{F}_1$ of the ellipsoid is parametrized by  
\[(A_1(u,v),B_1(u,v),C_1(u,v))\] where:

\begin{equation*}\aligned
A_1(u,v)=&\frac{\cos^{3}u \left( {a}^{2}-{b
}^{2} \right) }{a\sqrt{a^2+b^2}} \sqrt{  
 \left( {a}^{2}-{b}^{2} \right)\cos^{2}v + \left( {a}^{2}+{c}^{2}\right) 
  \sin ^{2}v }\\
B_1(u,v)=&-\frac {\sin ^{3}u \sin
 v \left(  {a}^{2}-{b}^{2} \right) }{b}\\
 C_1(u,v)=& \frac{\cos v }{c\sqrt{{a}^{2}+{c}^{2}}}  \left[  \left( {b}^{2}+{c}^{2} \right)  \cos^2 u + \left( {a}^{2}+{c}^{2} \right) 
  \sin^{2}u\right]^{\frac{3}{2}}
\endaligned\end{equation*}

The focal set $\mathcal{F}_2$ of the ellipsoid is parametrized by $(A_2(u,v),B_2(u,v),C_2(u,v))$ where:

\begin{equation*}\aligned
A_2(u,v)=&\frac{\cos u }{a\sqrt{a^2 + c^2}} \left[(a^2 - b^2)\cos^2v + (a^2 + c^2)\sin ^2v\right]^{\frac{3}{2}}\\
B_2(u,v)=& \frac { \sin
 u \sin ^{3}v\left(  {b}^{2}+{c}^{2} \right) }{b}\\
 C_2(u,v)=&  \frac{(b^2 + c^2)\cos ^3v }{c\sqrt{a^2 + c^2}}  \sqrt{(b^2 + c^2)cos^2 u + (a^2 + c^2)\sin^2u} 
\endaligned\end{equation*}
\end{pro}

\begin{figure}[H]
\centering
\includegraphics[width=0.8\textwidth,angle=0]{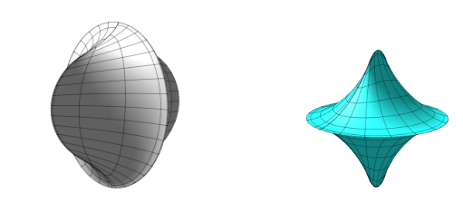}

\caption{ The focal surfaces $\mathcal{F}_1$ and $\mathcal{F}_2$ of the ellipsoid. Both are singular on two arcs of ellipses connecting the umbilic points and each is singular in an ellipse contained in a coordinate plane. At the umbilic points the singularities are of  type $D_4^{+}$ (Arnold's notation).}
\end{figure}

\begin{proof}
It follows directly from the parametrization of the ellipsoid $Z_{w}$
given by Remark \ref{rem:ellipsoid}.
    It is worth to observe that at the tropics defined by $\cos v=\pm c/\sqrt{b^2+c^2}$ the principal curvatures $k_i$ are unbounded but at these sets the normal $N$ has norm zero and the product $(1/k_i)N$ has a finite limit. See \cite{tari2012}.
\end{proof}

\bibliographystyle{plain}



\vspace{20mm}
\setlength{\parindent}{0pt}
\begin{minipage}{.50\linewidth}
\footnotesize
{\sc
Dimas N. T. Tejada\\
Facultad de Ciencias Naturales\\ y Matemática \\
  Universidad de El Salvador  \\
San Salvador, El Salvador \\
\textit{E-mail:} {\tt dimas.tejada@ues.edu.sv }\\
}
\end{minipage}
%
\begin{minipage}{.50\linewidth}
\footnotesize
{\sc
Ronaldo A. Garcia \\
Instituto de Matem\'atica e Estat\'istica \\
Universidade Federal de Goi\'as \\
Goiânia-GO, Brazil \\
\textit{E-mail:} {\tt ragarcia@ufg.br}
}
\end{minipage}

\end{document}